\documentclass[12pt]{amsart}

\usepackage{verbatim}
\usepackage{epsfig}

\newtheorem{thm}{Theorem}

\newtheorem*{claim}{Claim}
\newtheorem{conj}{Conjecture}
\newtheorem{prop}{Proposition}
\theoremstyle{definition}
\newtheorem*{rem}{Remark}
\newtheorem{ex}{Example}

\providecommand{\NN}{\mathbb{N}}
\providecommand{\ZZ}{\mathbb{Z}}

\providecommand{\RR}{\mathbb{R}}
\providecommand{\CC}{\mathbb{C}}

\def\eps{\varepsilon}
\def\tr{\mathrm{tr}}
\def\ext{\mathrm{ext}}

\begin{document}

\title[Lowest eigenvalues of differential operators]{A method for computing lowest eigenvalues of  symmetric 
polynomial differential operators by semidefinite programming}
\author{Jaka Cimpri\v c}
\date{submitted 25. 10. 2009, revised 15. 3. 2010}
\address{University of Ljubljana, Faculty of Mathematics and Physics, 
Department of Mathematics,  Jadranska 19, SI-1000 Ljubljana, Slovenija. 
e-mail: cimpric@fmf.uni-lj.si. www: http://www.fmf.uni-lj.si/\~{}cimpric/}

\begin{abstract}
A method for computing global minima of real multivariate polynomials based on
semidefinite programming was developed by N. Z. Shor, J. B. Lasserre and P. A. Parrilo.
The aim of this article is to extend a variant of their method to noncommutative symmetric polynomials 
in variables $X$ and $Y$ satisfying $YX-XY=1$ and $X^\ast=X$, $Y^\ast=-Y$.
Global minima of such polynomials are defined and showed to be equal to minima of the spectra
of the corresponding differential operators. We also discuss how to exploit sparsity
and symmetry. Several numerical experiments are included. The last section explains how our
theory fits into the framework of noncommutative real algebraic geometry.
\end{abstract}



\maketitle

\section{Motivation}
\label{sec1}

One of the most popular methods for computing global infima of real multivariate polynomials is
the method of sums of squares relaxations. The idea goes back to N. Z. Shor 
(his original papers are summarized in \cite[Chapter 9]{shor}) and it was further developed by J. B. Lasserre \cite{las} 
and P. A. Parrilo \cite{par}. See \cite{lau} for an extensive survey and \cite{sost,glob} for implementations. 
In this section we will present a variant of this method and in the next section we will extend it
from polynomials to polynomial differential operators. Later sections are concerned with  improvements of the basic method
and numerical experiments.

For a given polynomial $f \in \RR[\underline{X}]$,  $\underline{X}= (X_1,\ldots,X_d)$, write
\begin{equation}
\label{definf}
\inf f := \inf\{f(x) \mid x \in \RR^d\}=\sup \{ \mu \in \RR \mid f-\mu \ge 0 \mbox{ on } \RR^d\}
\end{equation}
for its unconstrained global infimum. By convention $\inf f=-\infty$ if $f$ is not bounded from below.
Note that the polynomial $X_1^2+(1-X_1 X_2)^2$ is bounded from below but it does not attain its infimum.

Write $\sum \RR[\underline{X}]^2$ for the set of all sums of squares of polynomials from $\RR[\underline{X}]$. 
Clearly, every element from $\sum \RR[\underline{X}]^2$ is nonnegative on $\RR^d$. In particular this is true for
the polynomial $s:=1+X_1^2+\ldots+X_d^2$. On the other hand, the Motzkin polynomial $1+X_1^2 X_2^4+X_1^4 X_2^2- 3 X_1^2 X_2^2$ 
is nonnegative on $\RR^d$, $d \ge 2$, but it does not belong to $\sum \RR[\underline{X}]^2$. 

Let us consider the following sequence of approximations of $\inf f$:
\begin{equation}
\label{defmu}
\mu_k(f) := \sup \{ \mu \in \RR \mid s^k (f-\mu) \in  \sum \RR[\underline{X}]^2\}, 
\end{equation}
where $k \in \NN:=\{0,1,2,\ldots\}$.
Since $\sum \RR[\underline{X}]^2$ is closed in the finest locally convex topology of 
$\RR[\underline{X}]$ by \cite[Corollary 11.6.4]{schbook}, we can replace $\sup$  by $\max$ in (\ref{defmu}).

\begin{prop} 
\label{prop1}
If $f \in \RR[\underline{X}]$ satisfies the following assumption:

\medskip

(*) \begin{minipage}{11cm}
$f_n(x) >0$ for every nonzero $x \in \RR^d$ where $n=2m$ is the total degree of $f$
and $f_n$ is the $n$-th homogeneous part of $f$.
\end{minipage}

\medskip

\noindent then $\lim \mu_k(f)=\inf f$.
\end{prop}

\begin{rem}
The assumption (*) is sufficent for the existence of global infimum but it is not necessary.
(A necessary condition is that $n=2m$ and $f_n(x) \ge 0$ for every nonzero $x \in \RR^d$.)
\end{rem}

\begin{proof}
It is clear from (\ref{definf}) and (\ref{defmu}) that $\inf f \ge \mu_k(f)$ for every $k$.
Since $s^k(f -\mu) \in \sum \RR[\underline{X}]^2$ implies
$s^{k+1}(f -\mu) \in \sum \RR[\underline{X}]^2$, we have by (\ref{defmu}) that
$\mu_{k+1}(f) \ge \mu_k(f)$ for every $k$. Finally, for every $\eps > 0$,
$f-\inf f+\eps$ is strictly positive on $\RR^d$. By Theorem \ref{thm1} 
we can find $k_\eps \in \NN$ such that 
$s^{k_\eps}(f-\inf f+\eps) \in \sum \RR[\underline{X}]^2$.
Hence, $\mu_{k_\eps}(f) \ge \inf f-\eps$ by (\ref{defmu}).
\end{proof}

The following dehomogenized version of a theorem of Reznick \cite{rez} was used in the proof,
see the comments after Theorem 5.5.2 in \cite{mm}.

\begin{thm}[Reznick 1995]
\label{thm1}
If $f \in \RR[\underline{X}]$ satisfies $f(x)>0$ for every $x \in \RR^d$ and the property (*)
then $s^k f \in \sum \RR[\underline{X}]^2$ for some $k \in \NN$.
\end{thm}

Finally, we would like to convince the reader that the numbers $\mu_k(f)$ can be effectively computed.
We will do so by reformulating the definition of $\mu_k(f)$ as a \textit{semidefinite program}, i.e. an optimization problem:
\begin{equation}
\label{defsdp}
\begin{minipage}{0.75 \textwidth}
minimize $\tr(CZ)$ subject to $\tr(A_iZ)=b_i$ and $Z\ge 0$
\end{minipage}
\end{equation}
where $C,A_i$ are given real symmetric matrices of the same size,
$b_i$ are given real numbers and $Z$ is an unknown real
symmetric matrix of the same size as $C,A_i$.

Let $v_k$ be a vector of all monomials of total degree less or equal to
$\frac12 \deg(s^k f) = k+m$. Its size is $n_k=\binom{k+m+d}{d}$.
Every element of $\sum \RR[\underline{X}]^2$ of degree  $\le 2(k+m)$
(in particular $s^k (f-\mu)$ for every $\mu \le \mu_k(f)$)
can be written in the form $v_k^T Z v_k$ where 
$Z$ is a positive semidefinite real symmetric matrix of size $n_k$. Therefore, by (\ref{defmu}),
 we can express $\mu_k(f)$ as the solution of the following optimization problem:

\begin{equation}
\label{musdp}
\begin{minipage}{0.75 \textwidth}
Find $\mu_k(f)=\max_{\mathcal F} \mu$ where $\mathcal{F}=\{(\mu,Z) \mid \mu \in \RR, $ \\
$Z \in M_{n_k}(\RR), \; Z^T=Z, \; Z\ge 0, \;  s^k (f-\mu)=v_k^T Z v_k\}$.
\end{minipage}
\end{equation}

Note that (\ref{musdp}) is not exactly a semidefinite program as defined by (\ref{defsdp}) but it can easily be converted into one by:
\begin{itemize}
\item eliminating $\mu$ using the linear relation obtained by comparing constant terms in $s^k (f-\mu)=v_k^T Z v_k$,
(we get $f_0-\mu=Z_{11}$ if $f_0$ is the constant term of $f$ and the first component of $v_k$ is $1$),
\item writing the linear relations between $Z_{ij}$ obtained by comparing coefficients in
$s^k (f-f_0+Z_{11})=v_k^T Z v_k$ in the form $\tr(A_i Z)=b_i$,
\item replacing $\mu_k(f)=\max \mu$ by $f_0-\mu_k(f)=\min Z_{11}=\min \tr(CZ)$ where $C_{11}=1$
and other $C_{ij}$ are $0$.
\end{itemize}

\begin{rem}
Recall that the Newton polytope $N(f)$ of a polynomial $f=\sum c_\alpha X^\alpha \in \RR[X]$
is the convex hull of its support $\{\alpha \mid c_\alpha \ne 0\}$. The main property of Newton polytopes
is $N(fg)=N(f)+N(g)$. (It is proved by showing that both sets have the same extreme points, i.e.
$\ext N(fg)=\ext N(f)+\ext N(g)=\ext (N(f)+N(g))$.) If $f=\sum g_i^2$, the property implies that
$N(g_i) \subseteq \frac12 N(f)$ for every $i$. Hence, the vectors $v_k$ from (\ref{musdp})
need not contain all monomials of degree $\le \frac12 \deg s^k f$ but only
the monomials from $\frac12 N(s^k(f-\mu))$.
\end{rem}

\section{Polynomial differential operators}

Our aim is to develop a similar theory for hermitian elements of the $d$-th Weyl algebra $\mathcal{W}(d)$.
Recall that $\mathcal{W}(d)$ is the unital complex $\ast$-algebra with generators $X_k,Y_l$,
defining relations $Y_l X_k-X_k Y_l=\delta_{kl}$ and involution $X_k^\ast=X_k$, $Y_l^\ast=-Y_l$,
where $k,l=1,\ldots d$. We will write $\mathcal{W}(d)_h$ for the
set of all elements $u \in \mathcal{W}(d)$ such that $u^\ast=u$ and
$\sum \mathcal{W}(d)^2$ for the set of all finite sums of elements $u^\ast u$ where $u \in \mathcal{W}(d)$.

The Schr\" odinger representation $\pi_0$ is the $\ast$-representation
of $\mathcal{W}(d)$ acting on the Schwartz space $\mathcal{S}(\RR^d)$ 
considered as dense domain of  $L^2(\RR^d)$,
which is defined by $(\pi_0(X_k) \phi)(t)=t_k \phi(t)$
and $(\pi_0(Y_l) \phi)(t)=\frac{\partial \phi}{\partial t_l}(t)$
for $k,l=1,\ldots,d$. We will write $\mathcal{W}(d)_+$ for the set of all elements $u \in \mathcal{W}(d)_h$ such that 
$\langle \pi_0(u) \phi, \phi \rangle \ge 0$ for all $\phi \in \mathcal{S}(\RR^d)$.
Clearly, $\sum \mathcal{W}(d)^2 \subseteq \mathcal{W}(d)_+$ while the converse is false
by \cite[Section 6]{sch}.

For a given element $c \in \mathcal{W}(d)_h$ write
\begin{equation}
\label{infc}
\inf c  :=  \sup\{ \lambda \in \RR \mid c-\lambda \cdot 1 \in \mathcal{W}(d)_+\} 
\end{equation}
if $\pi_0(c)$ is bounded from below. Otherwise write $\inf c=-\infty$. Clearly,
$\inf c$ is equal to the infimum of the numerical range of $\pi_0(c)$, i.e.
\[
\inf c =\inf \{\langle \pi_0(c)v,v \rangle \mid v \in \mathcal{S}(\RR^d),\Vert v\Vert=1\}.
\]
When $\pi_0(c)$ is bounded from below, one can define the Friedrichs extension $\pi_0(c)_F$
of $\pi_0(c)$ and show that 
\[
\inf c=\min \sigma(\pi_0(c)_F),
\]
see e.g. \cite{zettl} and the references therein.
If $\sigma(\pi_0(c)_F)$ is also discrete then $\inf c$ is equal to the lowest eigenvalue of $\pi_0(c)_F$.

\begin{rem}
There are well-known sufficient conditions on $V(X)$ implying that 
the Schr\" odinger operator $L=\pi_0(-\sum_{i=1}^d Y_i^2+V(X))$ is bounded from
below, is essentially selfadjoint (i.e. $L_F$ is the only selfadjoint extension of $L$)
and has discrete spectrum, see \cite[Sections 8.1, 8.2]{pan}.
\end{rem}

Let $c$ be an element of $\mathcal{W}(d)_h$ of even total degree $\deg c$.
We propose the following method for computing $\inf c$. Firstly, pick a sequence $s=(s_k)_{k \in \NN}$ of elements
of $\mathcal{W}(d)$ with $s_0=1$. Secondly, solve the following sequence of ``semidefinite programs'':
\begin{equation}
\label{muc}
\begin{minipage}{0.75 \textwidth}
Find $\mu_k(c,s)=\max_{\mathcal{F}} \mu$ where $\mathcal{F}=\{(\mu,Z) \mid \mu \in \RR, $ \\
$Z \in M_{n_k}(\CC), \; Z^H=Z, \; Z \ge 0, \; b_k^\ast (c-\mu) b_k=v_k^H Z v_k\}$
\end{minipage}
\end{equation}
where for every $k \in \NN$, $b_k=\prod_{i=0}^k s_k$, $v_k$ is a vector of all monomials $m$
in the generators of $\mathcal{W}(d)$ such that $\deg m \le \frac12 \deg b_k^\ast c b_k=\deg b_k+\frac12 \deg c$
and $n_k$ is the size of $v_k$.

Clearly, $\mu_{k+1}(c,s) \ge \mu_k(c,s)$ for every $k$ since $b_k$ divides $b_{k+1}$ in (\ref{muc}). 
If $\pi_0(s_k)$ is invertible for every $k$, then, by (\ref{infc}) and (\ref{muc}),
the sequence $\mu_k(c,s)$ is bounded above by $\inf c$.  The main question is what additional 
assumptions on $c$ and $s$ are needed to ensure that $\lim \mu_k(c,s)=\inf c$. Our numerical 
experiments suggest that the only additional assumption required is that $s_k$ are nonconstant
but we are unable to prove that. What we can prove about convergence is summarized in
Propositions 2 and 3 below; see also Conjecture \ref{conj1}.

Recall that \textit{the leading symbol} of an element 
\[
c = \sum_{\alpha,\beta \in \NN^d} c_{\alpha,\beta} X^\alpha Y^\beta
= \sum_{k=0}^{\deg c} \sum_{\alpha,\beta \in \NN^d \atop \vert \alpha \vert+\vert \beta \vert =k} c_{\alpha,\beta} X^\alpha Y^\beta
\in \mathcal{W}(d)
\]
(in multiindex notation) is the element
\[
\bar{c} = \sum_{\alpha,\beta \in \NN^d \atop \vert \alpha \vert+\vert \beta \vert =\deg c} c_{\alpha,\beta} X^\alpha \xi^\beta
\in \CC[X,\xi].
\]
If $c=c^\ast$, then $\bar{c}(X,i \xi) \in \RR[X,\xi]$. For example, the leading symbol of
$
N:= \sum_{i=1}^d \frac12 \left(X_i^2-Y_i^2-1 \right) \in \mathcal{W}(d)_h
$
is 
$
\bar{N}= \sum_{i=1}^d \frac12 \left(X_i^2-\xi_i^2\right).
$
Note that
$
\bar{N}(X,i \xi)=\sum_{i=1}^d \frac12 \left(X_i^2+\xi_i^2\right)
$
is $>0$ if $(X,\xi) \ne (0,0)$.

\begin{prop}
\label{prop2}
Suppose that $c \in \mathcal{W}(d)_h$, $4 | \deg c$ and $\bar{c}(X,i\xi) >0$
for every $X,\xi \in \RR^d$ with $(X,\xi) \ne (0,0)$. Pick $\alpha \in \RR^+\setminus \NN$ and a sequence
$m_k \in \ZZ$ in which every integer appears infinitely many times and write 
$s_0=1$, $s_k=N+(m_k+\alpha) \cdot 1$ for $k \ge 1$. 
Then $\lim \mu_k(c,s)=\inf c$.
\end{prop}

This is an immediate consequence of the fact that $\sigma(N)=\NN$ and the following result of Schm\" udgen
which is a noncommutative analogue of Theorem \ref{thm1}. (Assumption (1) can be relaxed slightly, 
see \cite[Th. 1.2]{nahas}, but this does not help us here.)

\begin{thm} 
\label{thm2}
\cite[Th. 1.1]{sch}
Suppose $c \in \mathcal{W}(d)_h$, $\deg c=2m$, satisfies
\begin{enumerate}
\item There exists $\eps > 0$ such that $c-\eps \cdot 1 \in \mathcal{W}(d)_+$.
\item $\bar{c}(X,i\xi) >0$ for every $X,\xi \in \RR^d$ with $(X,\xi) \ne (0,0)$.
\end{enumerate}
Finally, fix $\alpha \in \RR^+\setminus \NN$ and write $\mathcal{N}$ for the set 
of all finite products of elements $N+(\alpha+n) \cdot 1$, where $n \in \ZZ$.

If m is even, then there exists $b \in \mathcal{N}$ such that 
$bcb \in \sum \mathcal{W}(d)^2$. If $m$ is odd, then there exists 
$b \in \mathcal{N}$ such that $\sum_{j=1}^d b (X_j+Y_j) c (X_j-Y_j) b \in \sum \mathcal{W}(d)^2$.
\end{thm}

Proposition \ref{prop3} is a variant of Proposition \ref{prop2}. We need some notation.
We assume that $d=1$ and write $q=X$ and $p=-iY$.
Every nonzero $c \in \mathcal{W}(1)$ can be uniquely expressed as
$$
c=\sum_{j=0}^{d_1} \sum_{l=0}^{d_2} \gamma_{jl} p^j q^l
= \sum_{n=0}^{d_2} f_n(p) q^n = \sum_{k=0}^{d_1} g_k(q)p^k,
$$
where $f_{d_2} \ne 0$ and $g_{d_1}  \ne 0$. In this case we say that $c$ has multidegree $(d_1,d_2)$.
We fix two nonzero reals $\alpha$ and $\beta$. Let $\mathcal{S}$ be the monoid generated by $s_1=p-\alpha i$,
$s_2=q-\beta i$ and $s_1^\ast$, $s_2^\ast$. It is an Ore set.

\begin{prop}
\label{prop3}
Suppose that $c \in \mathcal{W}(1)_h$ has multidegree $(2m_1,2m_2)$, where $m_1,m_2 \in \NN$,
$\gamma_{2m_1,2m_2} \ne 0$ and $f_{2m_2}$ and $g_{2m_1}$ are positive on the real line.
Since $\mathcal{S}$ is countable we can number its elements, say $(u_i)_{i \in \NN}$,
assuming $u_0=1$. Write $s_0=1$ and let, for every $k \ge 0$, $s_{k+1} \in \mathcal{S}$ be the 
common right multiple of $s_k$ and $u_{k+1}$ that exists by the Ore property. 
Then $\lim \mu_k(c,s)=\inf c$.
\end{prop}

This is an immediate consequence of the following result of Schm\" udgen:

\begin{thm}
\label{thm3}
\cite[Theorem 5]{sch2}.
Let $c$ be a nonzero hermitian element of $\mathcal{W}(1)$ of multidegree $(2m_1,2m_2)$, where $m_1,m_2 \in \NN$.
Suppose that:
\begin{enumerate}
\item There exists a bounded self-adjoint operator $T>0$ on $L^2(\RR)$ such that $\pi_0(c) \ge T$.
\item $\gamma_{2m_1,2m_2} \ne 0$ and $f_{2m_2}$ and $g_{2m_1}$ are positive on the real line.
\end{enumerate}
Then there exists an element $s \in \mathcal{S}$ such that $s^\ast c s \in \sum \mathcal{W}(1)^2$.
\end{thm}

Note that neither Proposition \ref{prop2} nor Proposition \ref{prop3} cover
Schr\" odinger operators $-Y^2+V(X)$ with polynomial potential $V(X)$
of degree $>2$. This case however fits into the following conjecture:

\begin{conj}
\label{conj1}
Suppose that $c \in \mathcal{W}(1)_h$ is bounded from below 
and $s_k=a X+i1$ for every $k \in \NN$. Then $2|\deg c$ and $\lim \mu_k(c,s)=\inf c$.
\end{conj}

This conjecture is true if the following claim from \cite{jp} is true:

\begin{claim}
Let $\mathcal{A}$ be the algebra obtained from $\mathcal{W}(1)$ 
with the addition of the generator $(aX+i1)^{-1}$ (for $a \in \RR$)
and the commutation relation $[p,(a X+i1)^{-1}]=ia(a X+i1)^{-2}$.
Then every positive element $c \in \mathcal{A}$ has a quadratic sum factorization 
$c=\sum d_k^*d_k$ for some finite set of elements $d_k \in \mathcal{A}$. 
\end{claim}

The proof in \cite{jp} seems to have a gap (where they use a result of Schm\" udgen).

\section{Implementation, Numerical Examples for $d=1$}
\label{impl}

The computation of $\mu_k(c,s)$ was implemented as follows. Firstly, the input for semidefinite programs 
(i.e. the matrices $C,A_i$ and the numbers $b_i$) was computed by 
Mathematica\textsuperscript{\textregistered} (Wolfram Research) in rational (i.e. exact)
arithmetics. Linear relations among $Z_{ij}$ had to be solved before they were converted 
into the form $\tr(A_iZ)=b_i$ to ensure that the matrices $A_i$ were linearly independent.
Secondly, the input data was exported to Matlab\textsuperscript{\textregistered} (Mathworks)
where it was solved by either SeDuMi \cite{jos} or 
SDPT3 \cite{sdpt3} semidefinite programming solver (through the 
Yalmip interface \cite{lof}) in floating point arithmetics.

The problems of the basic method are illustrated by the following toy example in $d=1$:

\begin{ex}
\label{ex1}
Write $c=(2 N+1)^2=(X^2-Y^2)^2$. Clearly, $\inf c=1$ but pretend we don't know that.
The element $c$ satisfies the assumptions of Proposition \ref{prop2}.
We fix $\alpha=\frac12$ and $s_k=N+\alpha \cdot 1$. We will compute
approximations $\mu_0(c,s),\ldots,\mu_3(c,s)$ of $\inf c$ using sedumi
and sdpt3 respectively. The results are in Table 1. 
The first approximation is very good but higher approximations are getting
worse while the theory says they should be getting better.

\begin{table}
\begin{tabular}{|l|l|l|l|l|}
\hline
k & $\mu_k(c,s)$ (sedumi) & $\mu_k(c,s)$ (sdpt3) & $n_k$ & $m_k$ \\
\hline
0 & 0.999999999760 & 0.999999993360 & 6 & 15\\
1 & 1.000044013288 & 0.999999989593 & 15 & 45 \\
2 & 1.112894824977 & 0.999999941936 & 28 & 91 \\
3 & 73.69340728792 & 24.89729311234 & 45 & 153 \\
\hline
\end{tabular}
\medskip
\caption{The table of Example 1. The numbers $n_k$ and $m_k$ refer to the size and number
of matrices $A_i$ in semidefinite programs used for computing $\mu_k(c,s)$.}
\end{table}

\end{ex}

The lesson that we learn is that semidefinite programs should be kept as small as possible.
The most natural way to do this is to exploit sparsity. (Later we will also discuss how to
exploit symmetry.) For every
\[
c = \sum c_{\alpha,\beta} X^\alpha Y^\beta \in \mathcal{W}(1)
\]
write $N'(c)$ for the convex hull of the set
\[
\bigcup_{\alpha,\beta \in \NN \atop c_{\alpha,\beta} \ne 0} 
\{(\alpha-k,\beta-k) \mid k=0,\ldots,\min(\alpha,\beta)\}.
\]
The point is that because of the relation $YX-XY=1$ we must replace
$\{(\alpha,\beta)\}$ by $\{(\alpha-k,\beta-k) \mid k=0,\ldots,\min(\alpha,\beta)\}$.
As outlined in Section \ref{sec1}, one can prove the property
$N'(fg)=N'(f)+N'(g)$ for every $f,g \in \mathcal{W}(1)$. The property implies that the vectors
$v_k$ from the definition of $\mu_k(c,s)$ need not contain
all monomials $X^\alpha Y^\beta$ with $\alpha+\beta \le \deg b_k^\ast f b_k$
but only those from $\frac12 N'(b_k^\ast (c-\mu)  b_k)$.
This method works particularly well if $b_k$ depend only on $X$.

\begin{ex} 
\label{semilog}
Let $\lambda_0(\beta)$ be the lowest eigenvalue of $\pi_0(c_\beta)$ where
\[
c_\beta=-Y^2+X^2+\beta X^4.
\]
The values of $\lambda_0(\beta)$ for various $\beta$ were computed in 
\cite[Table 1]{ban}  to 15 decimals. We refer to his values as ``exact''. 

Let $\mu_k(\beta)$ be the solution of the semidefinite program (\ref{muc}) for
$c=c_\beta$, $b_k=(i-X)^k$ and $v_k=(1,X,\ldots,X^{k+2},Y,X Y,\ldots,X^k Y)^T$.
For each $\beta \in \{0.0001, 1, 10000\}$ and $k=0,1,\ldots,14$, we will compute $\mu_k(\beta)$ by sdpt3.
Finally, for each $\beta$ we draw the semi-log plot of the sequence 
of relative errors of $\mu_k(\beta)$ with respect to $\lambda_0(\beta)$,
i.e. the plot of the sequence $k \mapsto \log_{10} \vert \frac{\lambda_0(\beta)-\mu_k(\beta)}{\lambda_0(\beta)} \vert$.
The results are presented in Figure 1. In theory these plots should decrease to $-\infty$.
By Figure 1, they decrease only during first 8-10 steps. Similar results are obtained for
$b_k=(1+X^2)^k$ and $v_k=(1,X,\ldots,X^{2k+2},Y,\ldots,X^{2k} Y)^T$, see Figure 2.

\begin{figure}
\includegraphics[height=6cm]{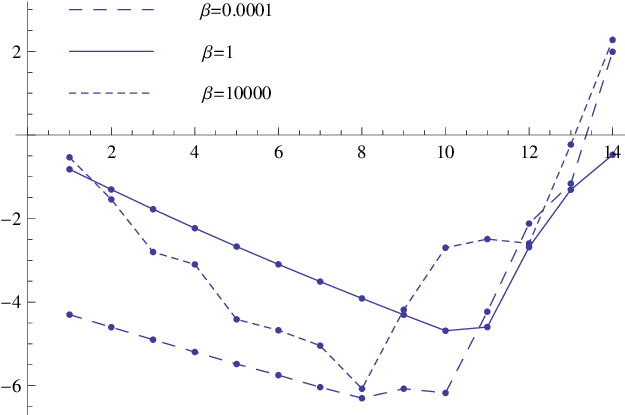}
\caption{The plots of $k \mapsto \log_{10} \vert \frac{\lambda_0(\beta)-\mu_k(\beta)}{\lambda_0(\beta)} \vert$ for $c=-Y^2+X^2+\beta X^4$,
$s_k=i-X$ and sparse $v_k$ (Example 2)}
\end{figure}
\begin{figure}
\includegraphics[height=6cm]{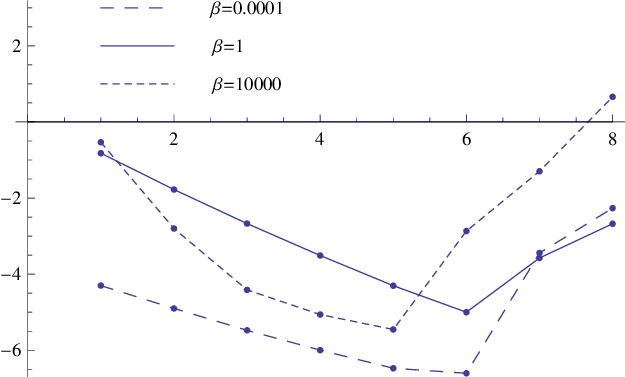}
\caption{The plots of $k \mapsto \log_{10} \vert \frac{\lambda_0(\beta)-\mu_k(\beta)}{\lambda_0(\beta)} \vert$ for $c=-Y^2+X^2+\beta X^4$,
$s_k=1+X^2$ and sparse $v_k$ (Example~2)}
\end{figure}
\end{ex}

Whatever improvement we use, things eventually start to go wrong (because of rounding errors)
and the question is when to stop.
If we use sdpt3, wrong values tend to undershoot, while if we use sedumi, they tend to overshoot.
We can use this observation to formulate an empirical stopping criterium for sdpt3: 
If $\mu_0 < \mu_1 < \ldots < \mu_l > \mu_{l+1}$, then return $\mu_l$ as the best approximation
for the lowest eigenvalue. We can also use $\min(\mu_l-\mu_{l-1},\mu_l-\mu_{l+1})$ as an estimate of its precision. 
There is no such stopping criterium for sedumi. We will use only sdpt3 in the sequel.

\begin{rem}
Another trick that sometimes improves stability in the commutative case is scaling
$X \to \lambda X$, $Y \to \lambda Y$. In our case this does not work, because it
violates the relation $Y X-X Y=1$. On the other hand, the transformation
$X \to \lambda X$, $Y \to \lambda^{-1} Y$ respects the relation but it does not 
improve stability.
\end{rem}

\section{Exploiting symmetry - finite groups}

Suppose that $G$ is a finite group acting on $\RR^d$
by orthogonal transformations. This action induces
in a natural way an action $\rho$ on the polynomial
ring $\RR[\underline{X}]$ and an action $\sigma$
on symmetric matrices that appear in our semidefinite
programs, see \cite{pg}. The action $\sigma$ helps
us put the matrices in our semidefinite programs into
block diagonal form and thus reduce the amount
of computation. It follows that every $G$-invariant
sum of squares is a sum of squares of invariant and 
semi-invariant polynomials, see \cite[Theorem 5.3]{pg}.

The same theory also works for Weyl algebras
and even some more general algebras, such us
enveloping algebras of Lie algebra.
Finite generation of the ring of invariants for
this situation was established in \cite{cck}
by passing to the corresponding graded ring.
An extension of \cite[Theorem 5.3]{pg} 
to locally finite-dimensions actions of compact groups
by $\ast$-automorphisms is provided by the last 
equation in the proof of \cite[Proposition 4]{ss}.

\begin{ex}
We would like to approximate the lowest eigenvalue of
\[
L=\pi_0(Y^4+X^4)
\]
by exploiting symmetry. Let $G=\{1,i,-1,-i\}$
act on $\mathcal{W}(1)$ by 
\[
\rho(i)(X)=iY \quad \text{and} \quad \rho(i)(Y)=iX.
\]
Clearly, $\rho(L)=L$. Since $\rho(i)(X+Y)=i(X+Y)$ and $\rho(i)(X-Y)=-i(X-Y)$,
it is more convenient to work with the generators
\[
a=\frac{X+Y}{\sqrt{2}} \quad \text{and} \quad a^\ast=\frac{X-Y}{\sqrt{2}}.
\]
We start with the zero-th approximation, i. e. we would like to find the
largest $\mu_0$ such that $Y^4+X^4-\mu_0 \cdot 1$ is a sum of hermitian squares.
We have to consider the 6-dimensional space $M_2$ of all monomials 
of degree less or equal to $2$. Eigenvectors of the restriction of $\rho$ to $M_2$ are:
\[\begin{array}{cl}
\lambda_1=1 : & 1,a^\ast a \\
\lambda_2=-1 : & a^2, (a^\ast)^2 \\
\lambda_3=i : & a,  \\
\lambda_4=-i : & a^\ast.
\end{array}\]
The $6 \times 6$ matrices that appear in our semidefinite program
can therefore be assumed to be block diagonal with two $2 \times 2$
blocks and two $1 \times 1$ blocks. We have therefore reduced 
the number of variables from $15$ to $8$. Using sdpt3, we get
\[
\mu_0=1.328 427 121
\]
To compute higher approximations, we need
apropriate denominators $b_k$ such that $b_k^\ast L b_k$ is still
$G$-invariant. Clearly, we can take for $b_k$ every $G$-semiinvariant
polynomial, e.g. any element from $\mathcal{N}$.
For 
\[
b_k=(2N+1)^k=(X^2-Y^2)^k=(2a^\ast a+1)^k
\]
we get
\[
\mu_1 = 1.396 727 721
\]
using the eigenvectors
\[\begin{array}{cl}
\lambda_1=1 : & 1, a^\ast a, a^4, (a^\ast)^4, (a^\ast)^2 a^2 \\
\lambda_2=-1 : & a^2, (a^\ast)^2, a^\ast a^3, (a^\ast)^3 a \\
\lambda_3=i : & a, (a^\ast)^3, a^\ast a^2  \\
\lambda_4=-i : & a^\ast, a^3, (a^\ast)^2 a \\
\end{array}\]
and 
\[
\mu_2=1.396 726 593
\]
using the eigenvectors
\[\begin{array}{cl}
\lambda_1=1 : & 1, a^\ast a, a^4, (a^\ast)^4, (a^\ast)^2 a^2, (a^\ast)^5 a, (a^\ast)^3 a^3, a^\ast a^5 \\
\lambda_2=-1 : & a^2, (a^\ast)^2, a^\ast a^3, (a^\ast)^3 a, a^6, (a^\ast)^2 a^4, (a^\ast)^4 a^2,(a^\ast)^6 \\
\lambda_3=i : & a, (a^\ast)^3, a^\ast a^2, (a^\ast)^4 a, (a^\ast)^2 a^3, a^5 \\
\lambda_4=-i : & a^\ast, a^3, (a^\ast)^2 a , a^\ast a^4, (a^\ast)^3 a^2, (a^\ast)^5 \\
\end{array}\]
Our stopping criterium tells us that $\mu_1$ is likely the best approximation we can get by this choice of $b_k$.
The method based on Conjecture \ref{conj1} gives
$\mu_4'=1.396 718 666$, $\mu_5'=1.396 726 819$, $\mu_6'=1.396 718 409$, i.e. a similar approximation
and a similar estimate for precision.
\end{ex}

\section{A conjecture about radial differential operators}

For every integer $d \ge 1$, we can identify the Hilbert space $L^2(\RR^+;r^{d-1})$ with the 
subspace of $L^2(\RR^d)$ consisting of radially invariant functions. Let $\mathcal{S}_d(\RR^+)$
be the subspace of $L^2(\RR^+;r^{d-1})$ which corresponds to the space of radially invariant
Schwartz functions on $\RR^d$. Let $\pi_d$ be the representation of the first Weyl algebra
$\mathcal{W}(1)$ which acts on $\mathcal{S}_d(\RR^+)$ by
$$(\pi_d(X)f)(r)=r f(r) \mbox{   and   } (\pi_d(Y)f)(r)=f'(r).$$
This representation is not a $\ast$-representation however if we consider the standard involution
on $\mathcal{W}(1)$ and the adjoint operation $L \mapsto L^+$ on differental operators. To make it
a $\ast$-representation we must consider a new involution $L \mapsto L^\ast := r^{d-1} L^+ r^{1-d}$ which
is conjugate to $L \mapsto L^+$.

The aim of this section is to provide numerical support for the following conjecture:

\begin{conj}
\label{conj2}
If $c \in \mathcal{W}(1)_h$ is such that $\pi_d(c) >0$ for some $d$ 
then $\deg_Y c$ is even and there exist $k \in \NN$  
and finitely many  $g_i,h_j \in \mathcal{W}(1)$ such that 
\[
(1+X^2)^k c \,(1+X^2)^k=\sum_i g_i^\ast g_i+ \sum_j h_j^\ast X h_j.
\]
\end{conj}

Suppose we want to compute the lowest eigenvalue $\lambda_0$ of a
$d$-dimensional radially invariant Schr\" odinger operator
\[
L=-\frac{1}{r^{d-1}} \frac{d}{dr} (r^{d-1} \frac{d}{dr}) + V(r)
\]
where $V(r)$ is a real polynomial in $r$ and $1/r$.
Pick the smallest $m \in \NN$ such that $r^{2m+d-1} V(r)$ has no negative powers.
Conjecture \ref{conj2} implies that the following sequence 
converges to $\lambda_0$:
\begin{equation}
\begin{minipage}{0.92 \textwidth}
$\begin{array}{r}
\mu_k  :=  \sup \{ \mu \in \RR \mid \exists g_i, h_j \in \mathcal{W}(1) :
(1+X^2)^k (-X^m Y X^{d-1} Y X^m + \\
 +X^{2m+d-1}(V(X)-\mu))(1+X^2)^k =\sum_i g_i^\ast g_i+ \sum_j h_j^\ast X h_j\}.
\end{array}$
\end{minipage}
\end{equation}
We can write $g_i^\ast g_i+ \sum_j h_j^\ast X h_j=u_k^\ast A u_k+v_k^\ast (X B) v_k$,
where $u_k,v_k$ are suitable vectors of monomials and $A,B$ are positive semidefinite 
complex hermitian matrices. Therefore, we can rewrite the definition of $\mu_k$ as a semidefinite program.

\begin{ex}
\label{exr1}
For $d=1,\ldots,10$ and
\[
V(r)=r+r^2+r^3,
\]
we will compute approximations $\mu_0,\ldots,\mu_5$ of $\lambda_0$ using 
\[
u_k= \left\{ \begin{array}{ll} (1,X,\ldots,X^{2k+1+\frac{d-1}{2}},Y,\ldots,X^{2k+\frac{d-1}{2}} Y)^T & d \mbox{ odd}, \\
(1,X,\ldots,X^{2k+1+\frac{d}{2}})^T & d \mbox{ even}. \end{array} \right.
\]
and
\[
v_k=\left\{ \begin{array}{ll}
(1,X,\ldots,X^{2k+1+\frac{d-1}{2}})^T & d \mbox{ odd}, \\
(1,X,\ldots,X^{2k+\frac{d}{2}},Y,\ldots,X^{2k-1+\frac{d}{2}}Y)^T& d \mbox{ even}. \end{array} \right.
\]
and compare them with \cite[Table 1]{hs}. The results are in Table 2. 

\begin{table}
\begin{tabular}{|l|l|l|l|l|l|l|l|}
\hline
$d$ & exact & $\mu_0$  & $\mu_1$ & $\mu_2$ & $\mu_3$ & $\mu_4$ & $\mu_5$ \\
\hline
2 & 3.5644 & 3.4973 & 3.5623 & 3.5643 & \underline{3.5644} & 3.5643 & 3.5630 \\
3 & 5.3066 & 5.2277 & 5.3046 & \underline{5.3065} & 5.3064 & 5.3058 & 5.3034\\
4 & 7.0746 & 7.0073 & 7.0730 & \underline{7.0746} & 7.0746 & 7.0677 & 7.0405 \\
5 & 8.8720 & 8.8187 & 8.8709 & \underline{8.8720} & 8.8720 & 8.8717 & 8.8713\\
6 & 10.6987 & 10.6488 & 10.6978 & \underline{10.6986} & 10.6986 & 10.6874 & 10.6030 \\
7 & 12.5534 & 12.5337 & 12.5523 & 12.5534 & \underline{12.5534} & 12.5533 & 12.5530 \\
8 & 14.4348 & 14.4229 & 14.4345 & \underline{14.4348} & 14.4342 & 14.3901 & 14.0928 \\
9 & 16.3415 & 16.3338 & 16.3413 & \underline{16.3414} & 16.3414 & 16.3414 & 16.3411 \\
10 & 18.2720 & 18.2664 & 18.2719 & \underline{18.2720} & 18.2712 & 18.2041 & 16.5905 \\
\hline
\end{tabular}
\medskip
\caption{The table of Example \ref{exr1}. For each $d$, the best approximation is underlined.
Exact values for $\lambda_0$ are from \cite[Table 1]{hs}.}
\end{table}

There is a problem with $d=1$. Namely, by  \cite[Table 1]{hs}, $\lambda_0$
is approximately $1.8306$, while we can show (in exact aritmetics)
that $\mu_0$ is approximately $1.9051$. A possible explanation is
that the eigen\-function corresponding to 1.8306 is not in $L^2(\RR^+)$.
This has nothing to do with Conjecture \ref{conj2}.

\end{ex}

\begin{ex} 
\label{exr2}
Suppose that $d=1$ and
\[
V(r)=-\frac{\lambda}{r}+r.
\]
We will compute approximations $\mu_0,\ldots,\mu_5$ of the lowest eigenvalue
using $u_k=(1,\ldots,X^{2k+1},Y,\ldots,X^{2k+1} Y)^T$
and $v_k=(1,\ldots,X^{2k+1})^T$ and compare them with \cite[Table I]{kw}.
The results are in Table 3.

\begin{table}
\begin{tabular}{|l|l|l|l|l|l|l|l|}
\hline
$\lambda$ & exact & $\mu_0$  & $\mu_1$ & $\mu_2$ & $\mu_3$ & $\mu_4$ & $\mu_5$\\
\hline
0.0 & 2.3381 & 1.8899 & 2.3193 & 2.3368 & \underline{2.3380} & 2.3300 & 2.2971 \\
0.2 & 2.1673 & 1.7277 & 2.1490 & 2.1661 & \underline{2.1672} & 2.1582 & 2.1198 \\
0.4 & 1.9885 & 1.5583 & 1.9706 & 1.9874 & \underline{1.9884} & 1.9783 & 1.9372 \\
0.6 & 1.8011 & 1.3810 & 1.7838 & 1.8001 & \underline{1.8010} & 1.7898 & 1.7349 \\
0.8 & 1.6044 & 1.1951 & 1.5878 & 1.6035 & \underline{1.6043} & 1.5915 & 1.5447 \\
1.0 & 1.3979 & 1.0000 & 1.3820 & 1.3971 & \underline{1.3978} & 1.3832 & 1.3280\\
1.2 & 1.1808 & 0.7949 & 1.1657 & 1.1801 & \underline{1.1807} & 1.1641 & 1.0800 \\
1.4 & 0.9526 & 0.5790 & 0.9383 & 0.9520 & \underline{0.9525} & 0.9360 & 0.8712\\
1.6 & 0.7127 & 0.3516 & 0.6992 & 0.7121 & \underline{0.7124} & 0.6937 & 0.6241 \\
1.8 & 0.4603 & 0.1119 & 0.4476 & 0.4597 & \underline{0.4599} & 0.4374 & 0.3691 \\
\hline
\end{tabular}
\medskip
\caption{The table of Example \ref{exr2}. For each $\lambda$, the best approximation is underlined.
Exact values for $\lambda_0$ (rounded from 15 to 4 decimals) are from \cite[Table I]{kw}.}
\end{table}

\end{ex}

\begin{ex}
Suppose that $d=1$ and
\[
V(r)=a r^2 +\frac{b}{r^2}
\]
where $a>0$, $b > -\frac14$. We apply one step of our method 
and divide the result by $r$ on both sides. We get
\[
-\frac{d^2}{dr^2}+V(r)-2 \sqrt{a}\,(1+\sqrt{b+\frac14})=g^\ast g, 
\]
where $g=\frac{d}{dr}+\sqrt{a}\, r-(\frac12+\sqrt{b+\frac14})\, r^{-1}$, which implies the inequality
\[
\lambda_0 \ge 2 \sqrt{a}\,(1+\sqrt{b+\frac14}).
\]
Similarly, if $d=1$ and
\[
V(r)=\frac{a}{r} +\frac{b}{r^2}
\]
where $a<0$, $b > -\frac14$, then we get as above
\[
-\frac{d^2}{dr^2}+V(r)+\frac{a^2}{(1+2\sqrt{b+\frac14})^2}=h^\ast h, 
\]
where $h=\frac{d}{dr}-\frac{a}{1+2\sqrt{b+\frac14}}-(\frac12+\sqrt{b+\frac14})\, r^{-1}$, hence
\[
\lambda_0 \ge -\frac{a^2}{(1+2\sqrt{b+\frac14})^2}.
\]
It is shown in \cite{is} that both inequalities for $\lambda_0$ are in fact equalities
but this is not clear from our method. If $a \ge 0$ in the second case then $\lambda_0 \ge 0$ because
\[
-\frac{d^2}{dr^2}+V(r)=\frac{a}{r}+
(\frac{d}{dr}-\frac{\frac12+\sqrt{b+\frac14}}{r})^\ast(\frac{d}{dr}-\frac{\frac12+\sqrt{b+\frac14}}{r}).
\]
\end{ex}

\section{Other $\ast$-algebras}

The aim of this short section is to outline a possible extension of our theory from $\RR[\underline{X}]$
and $\mathcal{W}(d)$ to other $\ast$-algebras.

Let $A$ be a finitely generated real or complex unital $\ast$-algebra and $\mathcal{R}$ a family of equivalence classes of
irreducible (possibly unbounded) $\ast$-representations of $A$. We can consider the elements of $A$ as
``polynomials'' and elements of $\mathcal{R}$ as (evaluations in) ``points'', see \cite{sch1}.
For every element $c \in A$ such that $c^\ast=c$ we can define
\[
\inf c := \sup \{ \mu \in \RR \mid \pi(c-\mu \cdot 1) \ge 0 \mbox{ for every } \pi \in \mathcal{R} \}.
\]
Clearly, our method for computing $\inf c$ can be applied to $A$ if:
\begin{itemize}
\item the monomials in the generators  are linearly independent and
\item it satisfies an analogue of Theorem \ref{thm1}.
\end{itemize}
Examples of such algebras are: 
\begin{itemize}
\item algebras of matrix polynomials \cite{av},
\item enveloping algebras of finite dimensional real Lie algebras \cite{schenv},
\item algebras of trigonometric polynomials \cite{ns},
\item finitely generated free real algebras \cite{he}, \cite{klep1} (also \cite{pna}, \cite{klep2}).
\end{itemize}

\subsection*{Acknowledgement} 
I would like to thank  Igor Klep for doing a part of programming
and Konrad Schm\" udgen for his comments.

\end{document}